\documentclass[11pt]{amsart}
\usepackage[all]{xy}
\usepackage{amsmath}
\usepackage{amsfonts}
\usepackage{amssymb}
\usepackage{latexsym}
\usepackage{graphicx}
\usepackage{color}
\usepackage{amscd}
\usepackage{epsfig}
\usepackage{cite}

\newtheorem{problem}{Problem}[section]
\newtheorem{definition}[problem]{Definition}
\newtheorem{lemma}[problem]{Lemma}
\newtheorem{theorem}[problem]{Theorem}

\newtheorem{corollary}[problem]{Corollary}

\title{ Generic twisted $T$-adic exponential sums of polynomials}
\author{Chunlei Liu}
\address{Department of Mathematics, Shanghai Jiao Tong
University, Shanghai 200240, P.R. China, E-mail: clliu@sjtu.edu.cn}
\author{Chuanze Niu}\address{School of
Mathematical Sciences, Beijing Normal University, Beijing 100875,
P.R. China, E-mail: niuchuanze@mail.bnu.edu.cn}
\begin{document}
\maketitle
\begin{abstract}
The twisted $T$-adic exponential sum associated to a polynomial in
one variable is studied. An explicit arithmetic polygon is proved to
be the generic Newton polygon of the $C$-function of the twisted
T-adic exponential sum. It gives the generic Newton polygon of the
$L$-functions of twisted $p$-power order exponential sums.
\end{abstract}



\section{Introduction}

Let $W$ be a Witt ring scheme of Witt vectors, $\mathbb{F}_q$ the field
of characteristic $p$ with $q$ elements,
$\mathbb{Z}_{q}=W(\mathbb{F}_q)$, and
$\mathbb{Q}_q=\mathbb{Z}_q[\frac{1}{p}]$.

Let $\triangle\supsetneq\{0\}$ be an integral convex polytope
 in $\mathbb{R}^n$, and $I$ the set of vertices of
 $\triangle$ different from the origin.
Let
$$f(x)=\sum\limits_{v\in \triangle}(a_vx^v,0,0,\cdots)\in
W(\mathbb{F}_q[x_1^{\pm1},\cdots,x_n^{\pm1}])\text{ with }
\prod_{v\in I}a_v\neq0,$$ where $x^v=x_1^{v_1}\cdots x_n^{v_n}$ if
$v=(v_1,\cdots,v_n)\in\mathbb{Z}^n$.

Let $T$ be a variable. Let $\mu_{q-1}$ be the group of $q-1$-th roots of unity in $\mathbb{Z}_q$ and
$\chi=\omega^{-u}$ with
$u\in\mathbb{Z}^n/(q-1)$ a fixed multiplicative
character of $(\mathbb{F}_q^{\times})^n$ into $\mu_{q-1},$  where
$\omega: x\rightarrow \hat{x}$ is the Teichm\"{u}ller character.

\begin{definition}
The sum
$$S_{f,u}(k,T)=\sum_{x\in(\mathbb{F}_{q^k}^{\times})^{n}}\chi({\rm Norm}_{\mathbb{F}_{q^k}/\mathbb{F}_q}(x))(1+T)^{Tr_{\mathbb{Q}_{q^k}/\mathbb{Q}_p}(\hat{f}(\hat{x}))}\in\mathbb{Z}_q[[T]]$$
is called a twisted $T$-adic exponential sum. And the function
$$L_{f,u}(s,T)=\exp(\sum_{k=1}^{\infty}S_{f,u}(k,T)\frac{s^k}{k})\in 1+s\mathbb{Z}_q[[T]][[s]]$$
is called an $L$-function of a twisted T-adic exponential sum.
\end{definition}

\begin{definition}
The function
$$C_{f,u}(s,T)=\exp(\sum_{k=1}^{\infty}-(q^k-1)^{-1}S_{f,u}(k,T)\frac{s^k}{k}),$$ is called a
 $C$-function of a twisted $T$-adic exponential sums.
\end{definition}

The L-function and C-function determine each other :
$$L_{f,u}(s,T)=\prod_{i=0}^n C_{f,u}(q^is,T)^{(-1)^{n-i+1}{n\choose i}},$$
and $$C_{f,u}(s,T)= \prod_{j=0}^{\infty} L_{f,u}(q^js,T)^{(-1)^{n-1}{n-j+1\choose j}}.$$
By the last identity, one sees that
$$C_{f,u}(s,T)\in 1+s\mathbb{Z}_q[[T]][[s]].$$

The $T$-adic exponential sums were first introduced by Liu-Wan \cite{LWn}.
We view $L_{f,u}(s, T)$ and $C_{f,u}(s, T)$ as power
series in the single variable $s$ with coefficients in the $T$-adic
complete field $\mathbb{Q}_q((T))$. The $C$-function $C_{f,u}(s,T)$ was
shown to be $T$-adic entire in $s$ by Liu-Wan \cite{LWn} for $u=0$ and
Liu \cite{Liu2} for all $u.$

Let $\zeta_{p^m}$ be a primitive $p^m$-th root of unity, and
$\pi_m=\zeta_{p^m}-1.$ Then $L_{f,u}(s,\pi_m)$ is the $L$-function
of the exponential sums $S_{f,u}(k,\pi_m)$ studied by
Adolphson-Sperber \cite{AS1,AS2,AS3,AS4} for $m=1$ and by Liu-Wei
\cite{LW} and
 Liu \cite{Liu1} for $m\geq1.$

Let $C(\triangle)$ be the cone generated by $\triangle$. There is a degree
function $\deg$ on $C(\triangle)$ which is
$\mathbb{R}_{\geq0}$-linear and takes values $1$ on each
co-dimension $1$ face not containing $0$. For $a\not\in
C(\triangle)$, we define $\deg(a)=+\infty$. Write $C_u(\triangle)=C(\triangle)\cap (u+(q-1)\mathbb{Z}^n)$ and
$M_u(\Delta)=\frac{1}{q-1}C_u(\triangle).$ Let $b$ be the least positive integer such that $p^bu\equiv u(\mod q-1).$
Order elements of $\bigcup_{i=0}^{b-1}M_{p^iu}(\Delta)$ so that $\deg(x_1)\leq\deg(x_2)\leq\cdots.$

\begin{definition}
The infinite u-twisted Hodge polygon $H_{\Delta,u}^{\infty}$ of $\Delta$ is the convex
function on $[0,+\infty]$ with initial point 0 which is linear between consecutive integers and whose slopes are
$$\frac{\deg(x_{bi+1})+\deg(x_{bi+2})+\cdots+\deg(x_{b(i+1)})}{b},~~ i=0,1,\cdots.$$
\end{definition}

Write NP for the short of Newton polygon.
Liu \cite{Liu2} proved the following Hodge bound for the $C$-function $C_{f,u}(s,T).$
\begin{theorem} We have
$$T-\text{adic NP of }C_{f,u}(s,T)\geq\text{ord}_p(q)(p-1)H_{\Delta,u}^{\infty}.$$
\end{theorem}

In the rest of this paper, we assume that $\Delta=[0,d]$ and
$a=\log_pq.$

Fix $0\leq u\leq q-2.$ Write $u=u_0+u_1p+\cdots+u_{a-1}p^{a-1}$ with $0\leq u_i\leq p-1.$
Then we have
$$\frac{u}{q-1}=-(u_0+u_1p+\cdots),~~~u_i=u_{b+i}.$$
Write $p^iu=q_i(q-1)+s_i$ for $i\in \mathbb{N}$ with $0\leq
s_i<q-1,$ then $s_{a-l}=u_l+u_{l+1}p+\cdots+u_{a+l-1}p^{a-1}$ for
$0\leq l\leq a-1$ and $s_i=s_{b+i}.$

\begin{lemma}
The infinite u-twisted Hodge polygon $H_{[0,d],u}^{\infty}$ is the convex
function on $[0,+\infty]$ with initial point 0 which is linear between consecutive integers and whose slopes are
$$\frac{u_0+u_1+\cdots+u_{b-1}}{bd(p-1)}+\frac{l}{d},~~l=0,1,\cdots.$$
\end{lemma}
\proof One observes that $\frac{s_i}{q-1}+l$ with $0\leq i\leq b-1$ is just
a permutation of $x_{bl+1},x_{bl+2},\cdots,x_{b(l+1)}.$ The lemma
follows.
\endproof

\begin{definition}
For any $i,n\in\mathbb{N},$ we define
$$\delta_{\in}^{(i)}(n)=\left\{
                        \begin{array}{ll}
                          1, & \hbox{}pl+u_{b-i}\equiv n(\mod d)~~\text{for some}~~ l<d\{\frac{n}{d}\}; \\
                          0, & \hbox{}otherwise,
                        \end{array}
                      \right.
$$ where $\{\cdot\}$ is the fractional part of a real number. We also define $\delta_{\in}^{(i)}(0)=0.$
\end{definition}

\begin{definition}
The arithmetic polygon $P_{[0,d],u}$ is the convex function on $[0,+\infty]$
with initial point 0 which is linear between consecutive integers and whose slopes are
$$\omega(n)=\frac{1}{b}\sum\limits_{i=0}^{b-1}(\lceil\frac{(p-1)n+u_{b-i}}{d}\rceil-\delta_{\in}^{(i)}(n)),~~n\in\mathbb{N},$$
where $\lceil\cdot\rceil$ is the least integer equal or greater than
a real number.
\end{definition}
Write
$$\{x\}^{'}=1+x-\lceil x\rceil=\left\{
                        \begin{array}{ll}
                          \{x\}, & \hbox{if } \{x\}\neq0, \\
                          1, & \hbox{if } \{x\}=0,
                        \end{array}
                      \right.
$$ where $\{\cdot\}$ is the fractional part of a real number.
Define$$\varepsilon(u)=\min\{d\{\frac{u_i}{d}\}^{'}|1\leq i\leq b\}.$$

In this paper, we shall prove the following theorems.
\begin{theorem}\label{main}We have
$$P_{[0,d],u}\geq (p-1)H_{[0,d],u}^{\infty}.$$
Moreover, they coincide at the point $d.$
\end{theorem}\begin{theorem}\label{main1}Let $f(x)=\sum\limits_{i=0}^d(a_ix^i,0,0,\cdots)$, and $p>4d-\varepsilon(u).$ Then
$$T-\text{adic NP of }C_{f,u}(s,T)\geq\text{ord}_p(q)P_{[0,d],u}. $$
\end{theorem}

\begin{theorem}\label{main2}Let $f(x)=\sum\limits_{i=0}^d(a_ix^i,0,0,\cdots)$, and $p>4d-\varepsilon(u).$ Then there is a
non-zero polynomial $H_u(y)\in\mathbb{F}_q[y_i\mid
i=0,1,\cdots,d]$ such that
$$T-\text{adic NP of }C_{f,u}(s,T)=\text{ord}_p(q)P_{[0,d],u}$$ if and only if
 $H_u((a_i)_{i=0,1,\cdots,d})\neq0.$
\end{theorem}

\begin{theorem}\label{main3}
Let $f(x)=\sum\limits_{i=0}^d(a_ix^i,0,0,\cdots),$ $p>4d-\varepsilon(u)$ and $m\geq1.$ Then
$$\pi_m-\text{adic NP of }C_{f,u}(s,\pi_m)\geq\text{ord}_p(q)P_{[0,d],u}$$ with equality holding if and only if
$H_u((a_i)_{i=0,1,\cdots,d})\neq0.$
\end{theorem}

By a result of W. Li \cite{Lw}, we see if $p\nmid d,$ $L_{f,u}(s,\pi_m)$ is a polynomial of degree $p^{m-1}d.$ From the above
theorem, we shall deduce  the following.
\begin{theorem}\label{main4}
Let $f(x)=\sum\limits_{i=0}^d(a_ix^i,0,0,\cdots),$ $p>4d-\varepsilon(u)$ and $m\geq1.$ Then
$$\pi_m-\text{adic NP of }L_{f,u}(s,\pi_m)\geq\text{ord}_p(q)P_{[0,d],u}~~ on~~[0,p^{m-1}d]$$ with equality holding if and only if
$H_u((a_i)_{i=0,1,\cdots,d})\neq0.$
\end{theorem}
Note that $\varepsilon(0)=d$ and $P_{[0,d],0}$ is just the arithmetic polygon $p_{[0,d]}$
defined in Liu-Liu-Niu\cite{LLN}. So the above results are
generalizations of the corresponding results in
Liu-Liu-Niu\cite{LLN}. Initiated by a conjecture of Wan \cite{DWn},
the asymptotic behavior of generic Newton polygon of
$L_{f,u}(s,\pi_m)$ with $m=1$ was studied by Zhu
\cite{Zh1,Zh2,Zhu3}, Blache-F\'erard \cite{BF}, and
Blache-F\'erard-Zhu \cite{BFZ}.

\section{Arithmetic estimate}
In this section, let
$R_i$ be a finite subset of $\{1,2,\cdots,a\}\times(\frac{s_i}{q-1}+\mathbb{N})$ with cardinality $an$ for each $1\leq i\leq b,$
$\tau$ a permutation of $R=\bigcup_{i=1}^bR_i.$

Write $i(l)=i$ if $l\in R_i.$ We shall estimate
$$\sum_{i=1}^b\sum_{l\in
R_i}\lceil\frac{p\phi_i(l)-\phi_{i(\tau(l))}(\tau(l))+u_{b-i}}{d}\rceil,$$ where
$\phi_i$ is the projection
$\{1,2,\cdots,a\}\times(\frac{s_i}{q-1}+\mathbb{N})\rightarrow \mathbb{N}$ such that $$\phi_i((\cdot,\frac{s_{i}}{q-1}+l))=l.$$

Define $$\delta_{<}^{(i)}(l)=\left\{
                                \begin{array}{ll}
                                  1, & \hbox{} if~~\{\frac{l}{d}\}^{'}<\{\frac{pl+u_{b-i}}{d}\}^{'};\\
                                  0, & \hbox{}otherwise.
                                \end{array}
                              \right.
$$

Write $A_n=\{0,1,\cdots,n-1\}$ and $\mu=\{\frac{n-1}{d}\}.$
\begin{lemma}\label{mainrelation}
For any $1\leq i\leq b,$ we have
$$\sum_{l=0}^{n-1}(\delta_{<}^{(i)}(l)-\delta_{\in}^{(i)}(l))$$
$$=\sharp\{l\in A_n|\{\frac{l}{d}\}^{'}\leq
\mu<\{\frac{pl+u_{b-i}}{d}\}^{'}\} -\sharp\{l\in A_n|\{\frac{l}{d}\}^{'}>
\mu\geq\{\frac{pl+u_{b-i}}{d}\}^{'}\}.
$$
\end{lemma}
\proof Note that both $\delta_{<}^{(i)}$ and $\delta_{\in}^{(i)}$
have a period $d$ and initial value 0, so we may assume that $n\leq
d.$ The case $n=1$ is trivial. We are going to show this for
$n\geq2.$ Write
$$\delta_{[m,n]}^{(i)}=\left\{
                                      \begin{array}{ll}
                                        1, & \hbox{} pl+u_{b-i}\equiv0(\mod d)~~for~~some~~ m\leq l\leq n ;\\
                                        0, & \hbox{}otherwise.
                                      \end{array}
                                    \right.
$$
Hence by definition, we have
$$\sum_{l=0}^{n-1}\delta_{<}^{(i)}(l)=\sharp\{1\leq l\leq
n-1|\{\frac{l}{d}\}<\{\frac{pl+u_{b-i}}{d}\}\}+\delta_{[1,n-1]}^{(i)}.$$
We also have
$$\sum_{l=0}^{n-1}\delta_{\in}^{(i)}(l)=\sharp\{1\leq l\leq n-2|l+1\leq d\{\frac{pl+u_{b-i}}{d}\}\leq n-1\}+1_{\frac{1}{d}\leq\{\frac{u_{b-i}}{d}\}\leq\frac{n-1}{d}}.$$
Therefore
$$
\sum_{l=0}^{n-1}(\delta_{<}^{(i)}(l)-\delta_{\in}^{(i)}(l))$$
$$=\sharp\{1\leq l\leq
n-1|\{\frac{pl+u_{b-i}}{d}\}>\mu\}+\delta_{[1,n-1]}^{(i)}-1_{\frac{1}{d}\leq\{\frac{u_{b-i}}{d}\}\leq\frac{n-1}{d}}.$$
Note that $$\sharp\{l\in A_n|\{\frac{l}{d}\}^{'}\leq
\mu<\{\frac{pl+u_{b-i}}{d}\}^{'}\}$$
$$=\sharp\{1\leq l\leq
n-1|\{\frac{pl+u_{b-i}}{d}\}>\mu\}+\delta_{[1,n-1]}^{(i)}$$
and $$\sharp\{l\in A_n|\{\frac{l}{d}\}^{'}>
\mu\geq\{\frac{pl+u_{b-i}}{d}\}^{'}\}=1_{\frac{1}{d}\leq\{\frac{u_{b-i}}{d}\}\leq\frac{n-1}{d}},$$
the lemma follows.
\endproof

\begin{lemma} Let $A,B,C$ be sets
with $A$ finite, and $\tau$ a permutation of $A$. Then
$$\#\{a\in A\mid \tau(a)\in B,a\in C\}
$$$$\geq\#\{a\in A\mid a\in B,a\in C\}-\#\{a\in A\mid a\not\in
B,a\not\in C\}.$$
\end{lemma}
\proof The reader may refer \cite {LLN} and we omit the proof here.\endproof
 For $l\in R_i,$ define
$$\delta_{<,\tau}^{(i)}(l)=\left\{
                                \begin{array}{ll}
                                  1, & \hbox{} if~~\{\frac{\phi_{i(\tau(l))}(\tau(l))}{d}\}^{'}<\{\frac{p\phi_i(l)+u_{b-i}}{d}\}^{'};\\
                                  0, & \hbox{}otherwise.
                                \end{array}
                              \right.
$$
\begin{theorem}\label{estiofdelta}
We have
$$
\sum_{i=1}^{b}\sum_{l\in R_i}\delta_{<,\tau}^{(i)}(l)\geq
a\sum_{i=1}^b\sum_{l=0}^{n-1}(\delta_{<}^{(i)}-\delta_{\in}^{(i)})(l)-2\sum_{i=1}^b\sharp\{l\in
R_i|\phi_i(l)>n-1\}.
$$
\end{theorem}
\proof  By
definition, we have
$$\sum_{i=1}^{b}\sum_{l\in
R_i}\delta_{<,\tau}^{(i)}(l)=\sharp\{l\in
R|\{\frac{\phi_{i(\tau(l))}(\tau(l))}{d}\}^{'}<\{\frac{p\phi_{i(l)}(l)+u_{b-i(l)}}{d}\}^{'}\}
$$
$$\geq\sharp\{l\in
R|\{\frac{\phi_{i(\tau(l))}(\tau(l))}{d}\}^{'}\leq\mu<\{\frac{p\phi_{i(l)}(l)+u_{b-i(l)}}{d}\}^{'}\}.
$$
 Applying the last lemma with $A=R$ and
$$B=\{l\in R|\{\frac{\phi_{i(l)}(l)}{d}\}'\leq\mu\},~C=\{l\in
R|\mu<\{\frac{p\phi_{i(l)}(l)+u_{b-i(l)}}{d}\}'\},$$ we have
$$
\sum_{i=1}^{b}\sum_{l\in
R_i}\delta_{<,\tau}^{(i)}(l)\geq\sharp\{l\in
R|\{\frac{\phi_{i(l)}(l)}{d}\}^{'}\leq\mu<\{\frac{p\phi_{i(l)}(l)+u_{b-i(l)}}{d}\}^{'}\}$$
$$-\sharp\{l\in
R|\{\frac{\phi_{i(l)}(l)}{d}\}^{'}>\mu\geq\{\frac{p\phi_{i(l)}(l)+u_{b-i(l)}}{d}\}^{'}\}
$$ Note that
$$\sharp\{l\in
R_i|\{\frac{\phi_i(l)}{d}\}^{'}\leq\mu<\{\frac{p\phi_i(l)+u_{b-i}}{d}\}^{'}\}$$
$$\geq\sharp\{\phi_i(l)\in A_n|\{\frac{\phi_i(l)}{d}\}^{'}\leq\mu<\{\frac{p\phi_i(l)+u_{b-i}}{d}\}^{'}\}
-\sharp\{\phi_i(l)\in A_n,l\not\in R_i\}$$
$$=a\sharp\{l\in A_n|\{\frac{l}{d}\}^{'}\leq
\mu<\{\frac{pl+u_{b-i}}{d}\}^{'}\}-\sharp\{\phi_i(l)\in A_n,l\not\in
R_i\},
$$
and
$$\sharp\{l\in
R_i|\{\frac{\phi_i(l)}{d}\}^{'}>\mu\geq\{\frac{p\phi_i(l)+u_{b-i}}{d}\}^{'})$$
$$\leq\sharp\{\phi_i(l)\in A_n|\{\frac{\phi_i(l)}{d}\}^{'}>\mu\geq\{\frac{p\phi_i(l)+u_{b-i}}{d}\}^{'}\}+\sharp\{l\in
R_i,\phi_i(l)>n-1\}$$ $$ =a\sharp\{l\in A_n|\{\frac{l}{d}\}^{'}>
\mu\geq\{\frac{pl+u_{b-i}}{d}\}^{'}\}+\sharp\{l\in
R_i,\phi_i(l)>n-1\}.
$$
By $\sharp\{\phi_i(l)\in A_n,l\not\in R_i\}=\sharp\{l\in
R_i,\phi_i(l)\not\in A_n\},$ the theorem follows from Lemma \ref{mainrelation}.
\endproof

\begin{theorem}\label{mainesti}
If $p>4d-\varepsilon(u),$ then we have
$$\sum_{i=1}^b\sum_{l\in
R_i}\lceil\frac{p\phi_i(l)-\phi_{i(\tau(l))}(\tau(l))+u_{b-i}}{d}\rceil
 \geq  abP_{[0,d],u}(n)+\sharp\{l\in R,\phi_{i(l)}(l)>n-1\}.$$
\end{theorem}
\proof  By the definition of $\delta_{<,\tau}^{(i)}(l),$ we see
$$\lceil\frac{p\phi_i(l)-\phi_{i(\tau(l))}(\tau(l))+u_{b-i}}{d}\rceil=
\lceil\frac{p\phi_i(l)+u_{b-i}}{d}\rceil-\lceil\frac{\phi_{i(\tau(l))}(\tau(l))}{d}\rceil+\delta_{<,\tau}^{(i)}(l).$$
Since $p>4d-\varepsilon(u),$ we have $$\sum_{l\in
R_i}(\lceil\frac{p\phi_i(l)+u_{b-i}}{d}\rceil-\lceil\frac{\phi_i(l)}{d}\rceil)$$
$$=(\sum_{\phi_i(l)\in A_n}+\sum_{\stackrel{l\in
R_i}{\phi_i(l)\not\in A_n}}-\sum_{\stackrel{\phi_i(l)\in A_n}{l\not\in
R_i}})(\lceil\frac{p\phi_i(l)+u_{b-i}}{d}\rceil-\lceil\frac{\phi_i(l)}{d}\rceil)$$
$$\geq a\sum_{l=0}^{n-1}(\lceil\frac{pl+u_{b-i}}{d}\rceil-\lceil\frac{l}{d}\rceil)+3\sharp\{l\in
R_i|\phi_i(l)>n-1\}.
$$
The theorem is proved by Theorem \ref{estiofdelta}.
\endproof

\section{Twisted T-adic Dwork's trace formula}
In this section we review the twisted $T$-adic analogy of Dwork's theory on exponential sums.

Let
$E(x)=\exp(\sum\limits_{i=0}^{\infty}\frac{x^{p^i}}{p^i})=\sum\limits_{n=0}^{\infty}\lambda_nx^n$
be the Artin-Hasse series. Define a new T-adic uniformizer $\pi$ of
the $T$-adic local ring $\mathbb{Q}_p[[T]]$ by the formula
$E(\pi)=1+T.$

Recall $M_u=M_u([0,d])=\{v\in\mathbb{Q}_{\geq0}|v\equiv \frac{u}{q-1}(\mod 1)\}.$ Write
$$L_u=\{\sum_{v\in
M_u}b_v\pi^{\frac{v}{d}}x^{v}:b_v\in\mathbb{Z}_q[[\pi^{\frac{1}{d(q-1)}}]]\}$$ and
$$B_u=\{\sum_{v\in
M_u}b_v\pi^{\frac{v}{d}}x^{v}:b_v\in\mathbb{Z}_q[[\pi^{\frac{1}{d(q-1)}}]];{\rm
ord}_{\pi}b_v\rightarrow\infty,v\rightarrow\infty\}.$$

Define $\psi_p:L_u\rightarrow L_{p^{-1}u}$ by
$\psi_p(\sum\limits_{v\in M_u}b_vx^v)=\sum\limits_{v\in M_{p^{-1}u}}b_{pv}x^v.$ Then the map $\psi_p\circ E_f$
sends $L_u$ to $B_{p^{-1}u},$ where $E_f(x)=\prod\limits_{a_i\neq0}E(\pi\hat{a}_ix^i)\in L_0.$

The Galois group ${\rm Gal}(\mathbb{Q}_q/\mathbb{Q}_p)$ can act on
$B_u$ by fixing $\pi^{\frac{1}{d(q-1)}}$ and $x.$ Let $\sigma\in{\rm
Gal}(\mathbb{Q}_q/\mathbb{Q}_p)$ be the Frobenius element such that
$\sigma(\zeta)=\zeta^p$ if $\zeta$ is a $(q-1)$-th root of unity.
The operator $\Psi=\sigma^{-1}\circ\psi_p\circ E_f(x)$ sends $B_u$
to $B_{p^{-1}u},$ hence $\Psi$ operators on
$B=\bigoplus_{i=0}^{b-1}B_{p^iu}.$ We call it Dwork's
$T$-adic semi-linear operator because it is semi-linear over
$\mathbb{Z}_q[[\pi^{\frac{1}{d(q-1)}}]].$

Note that $\Psi^a=\psi_p^a\circ
\prod_{i=0}^{a-1}E_f^{\sigma^i}(x^{p^i}).$ It follows that $\Psi^a$
operates on $B_u$ and is linear over
$\mathbb{Z}_p[[\pi^{\frac{1}{d(q-1)}}]].$ Moreover, one can show
that $\Psi$ is completely continuous in sense of Serre \cite{Se}, so
$\det(1-\Psi^as|B_u/\mathbb{Z}_q[[\pi^{\frac{1}{d(q-1)}}]])$ and
$\det(1-\Psi s|B/\mathbb{Z}_p[[\pi^{\frac{1}{d(q-1)}}]])$ are well
defined.

Now we state the twisted $T$-adic Dwork's trace formula \cite{Liu1}.
\begin{theorem}We have
$$C_{f,u}(s,T)=\det(1-\Psi^as|B_u/\mathbb{Z}_q[[\pi^{\frac{1}{d(q-1)}}]]).$$
\end{theorem}

\section{The twisted Dwork semi-linear operator}
In this section, we shall study the twisted Dwork semi-linear operator
$\Psi.$

Recall that $\Psi=\sigma^{-1}\circ\psi_p\circ E_f,$
$E_f(x)=\prod_{i=0}^{d}E(\pi\hat{a}_ix^i)=\sum\limits_{n=0}^{\infty}\gamma_nx^n,$
where
$$\gamma_n=\sum\limits_{\sum\limits_{i=0}^{d}in_i=n;n_i\geq0}\pi^{\sum\limits_{i=0}^dn_i}
\prod\limits_{i=0}^d\lambda_{n_i}\hat{a}_i^{n_i}.$$

We see $B=\bigoplus_{i=1}^bB_{p^iu}$ has a basis represented by
$\{x^{\frac{s_i}{q-1}+j}\}_{1\leq i\leq b,~j\in\mathbb{N}}$ over $\mathbb{Z}_q[[\pi^{\frac{1}{d(q-1)}}]].$
 We have
$$\Psi(x^{\frac{s_i}{q-1}+j})=\sigma^{-1}\circ\psi_p(\sum_{l=0}^{\infty}\gamma_lx^{\frac{s_i}{q-1}+j+l})
=\sum_{l=0}^{\infty}\gamma_{pl+u_{b-i}-j}^{\sigma^{-1}}x^{l+\frac{s_{i-1}}{q-1}}.$$

For $1\leq k,i\leq b$ and $l,j\in\mathbb{N},$ define
$$\gamma_{(\frac{s_k}{q-1}+l,\frac{s_i}{q-1}+j)}=\left\{
                        \begin{array}{ll}
                          \gamma_{pl+u_{b-i}-j}, & \hbox{} k=i-1;\\
                          0, & \hbox{} otherwise.
                        \end{array}
                      \right.
$$
Let
$\xi_1,\cdots,\xi_{a}$ be a normal basis of $\mathbb{Q}_q$ over
$\mathbb{Q}_p$ and write
$$\xi_v^{\sigma^{-1}}\gamma_{(\frac{s_k}{q-1}+l,\frac{s_i}{q-1}+j)}^{\sigma^{-1}}=
\sum_{w=1}^{a}\gamma_{(w,\frac{s_k}{q-1}+l)(v,\frac{s_i}{q-1}+j)}\xi_w.$$
It is easy to see $\gamma_{(w,\frac{s_k}{q-1}+l)(v,\frac{s_i}{q-1}+j)}=0$ for any $w$ and $v$ if $k\neq i-1.$
 Define the $i$-th submatrix by
$$\Gamma^{(i)}=(\gamma_{(w,\frac{s_{i-1}}{q-1}+l)(v,\frac{s_i}{q-1}+j)})_{1\leq w,v\leq a;l,j\in\mathbb{N}},$$
then the matrix of the operator $\Psi$ on B over
$\mathbb{Z}_p[[\pi^{\frac{1}{d(q-1)}}]]$ with respect to the basis
$\{\xi_vx^{\frac{s_i}{q-1}+j}\}_{1\leq i\leq b,1\leq v\leq a;j\in\mathbb{N}}$ is

$$\Gamma=\left(
                \begin{array}{ccccc}
                  0 & \Gamma^{(1)} & 0 & \cdots & 0 \\
                  0 & 0 & \Gamma^{(2)} & \cdots & 0 \\
                  \vdots & \vdots & \vdots & \vdots & \vdots \\
                  0 & 0 & 0 & \cdots & \Gamma^{(b-1)} \\
                  \Gamma^{(b)}& 0 & 0 & \cdots & 0 \\
                \end{array}
              \right).
$$

Hence by a
result of Li-Zhu \cite{LZ}, we have $$\det(1-\Psi
s|B/\mathbb{Z}_p[[\pi^{\frac{1}{d(q-1)}}]])=\det(1-\Gamma s)=\sum_{n=0}^{\infty}(-1)^{bn}C_{bn}s^{bn},$$
with $C_n=\sum\det(A)$ where
$A$ runs over all principle $n\times n$ submatrix  of $\Gamma.$

For every principle submatrix $A$ of $\Gamma,$ write
$A^{(i)}=A\cap\Gamma^{(i)}$ as the sub-matrix of $\Gamma^{(i)}.$ For
a principle $bn\times bn$ submatrix $A$ of $\Gamma,$ by linear
algebra, if one of $A^{(i)}$ is not $n\times n$ submatrix of
$\Gamma^{(i)}$ then at least one row or column of $A$ are 0 since
$A$ is principle. Let $\mathcal{A}_n$ be the set of all $bn\times bn$ principle
submatrix A of $\Gamma$ with $A^{(i)}$ all $n\times n$ submatrix of
$\Gamma^{(i)}$ for each $1\leq i\leq b.$ Then we have
$$C_{bn}=\sum_{A\in\mathcal{A}_n}\det(A)=\sum_{A\in\mathcal{A}_n}(-1)^{n^2(b-1)}\prod_{i=1}^b\det(A^{(i)}).$$

Let $O(\pi^{\alpha})$ denotes any element of $\pi$-adic order $\geq
\alpha.$
\begin{lemma}\label{rawesti} We have
$$\gamma_n=\pi^{\lceil\frac{n}{d}\rceil}\sum\limits_{\stackrel{\sum\limits_{i=1}^din_i=n}
{\sum\limits_{i=1}^dn_i =\lceil\frac{n}{d}\rceil}}
\prod_{i=1}^d\lambda_{n_i}\hat{a}_i^{n_i}
+O(\pi^{\lceil\frac{n}{d}\rceil+1}).$$
\end{lemma}

\proof This follows from the fact that
$\sum\limits_{i=0}^dn_i\geq\lceil\frac{n}{d}\rceil$ if
$\sum\limits_{i=0}^din_i=n.$
\endproof

\begin{corollary}\label{ord of gamma}
For any $1\leq i\leq b$ and $1\leq w,v\leq a,$ we have
$$\gamma_{(w,\frac{s_{i-1}}{q-1}+l)(v,\frac{s_i}{q-1}+j)}
=O(\pi^{\lceil\frac{pl+u_{b-i}-j}{d}\rceil}).$$
\end{corollary}

\proof By
$$\xi_v^{\sigma^{-1}}\gamma_{(\frac{s_{i-1}}{q-1}+l,\frac{s_i}{q-1}+j)}^{\sigma^{-1}}=\xi_v^{\sigma^{-1}}\gamma_{pl-j+u_{b-i}}^{\sigma^{-1}}
=\sum_{w=1}^{a}\gamma_{(w,\frac{s_{i-1}}{q-1}+l)(v,\frac{s_i}{q-1}+j)}\xi_w,$$ this follows from
Lemma \ref{rawesti}.\endproof

\begin{theorem}\label{order of C_bn}
If $p>4d-\varepsilon(u),$ then we have $$\text{ord}_{\pi}(C_{abn})\geq
abP_{[0,d],u}(n).$$ In particular, $$C_{abn}=\pm{\rm
Norm}(\prod_{i=1}^b\det((\gamma_{pl-j+u_{b-i}})_{l,j\in
A_n}))+O(\pi^{abP_{[0,d],u}(n)+\frac{1}{d(q-1)}}),$$
where Norm is the norm map from
$\mathbb{Q}_q(\pi^{\frac{1}{d(q-1)}})$ to
$\mathbb{Q}_p(\pi^{\frac{1}{d(q-1)}})$.
\end{theorem}
\proof Let $A\in \mathcal{A}_{an},$ $R_i$ the set of rows of
$A^{(i)}$ as the submatrix of $\Gamma^{(i)},$ $\tau$ a permutation
of $R=\bigcup_iR_i.$ By the above corollary and Theorem
\ref{mainesti}, we have
$$\text{ord}_{\pi}(C_{abn})\geq\sum_{i=1}^b\sum_{l\in R_i}\lceil\frac{p\phi_i(l)+u_{b-i}-\phi_{i(\tau(l))}(\tau(l))}{d}\rceil
\geq abP_{[0,d],u}(n).$$ Moreover the strict inequality holds if
there exist $1\leq i\leq b$ such that $R_i\neq\phi_i^{-1}(A_n),$ hence
$$C_{abn}=\prod_{i=1}^b\det(\gamma_{(w,\frac{s_{i-1}}{q-1}+l)(v,\frac{s_i}{q-1}+j)})_{1\leq w,v\leq a;l,j\in A_n}+O(\pi^{abP_{[0,d],u}(n)+\frac{1}{d(q-1)}}).$$
Therefore the theorem follows from
$$\det((\gamma_{(w,\frac{s_{i-1}}{q-1}+l)(v,\frac{s_i}{q-1}+j)})_{1\leq w,v\leq a;l,j\in A_n}=\pm{\rm Norm}(\det(\gamma_{(\frac{s_{i-1}}{q-1}+l,\frac{s_i}{q-1}+j)})_{l,j\in A_n}).$$
\endproof

\section{Hasse polynomial}
In this section, $1\leq n\leq d-1.$ We shall study
$\det((\gamma_{pl-j+u_{b-i}})_{0\leq l,j\leq n-1})).$

\begin{definition}Let $S_n$ be the set of permutations of $A_n=\{0,1,\cdots,n-1\}.$
For $1\leq i\leq b,$ we define
$$S_{n,i}=\{\tau\in
S_n|\frac{\tau(l)}{d}\geq\frac{pl+u_{b-i}}{d}-\lceil\frac{pl+u_{b-i}-(n-1)}{d}\rceil\}.$$
\end{definition}
\begin{lemma}\label{estiofdet}
Let $p>4d-\varepsilon(u)$ and $\tau\in S_n.$   Then we have
$$\sum_{i=1}^b\sum_{l=0}^{n-1}\lceil\frac{pl+u_{b-i}-\tau(l)}{d}\rceil\geq
bP_{[0,d],u}(n),$$ with equality holding if and only if $\tau\in
S_{n,i}$ for each $1\leq i\leq b.$
\end{lemma}

\proof Because
$$\lceil\frac{pl+u_{b-i}-\tau(l)}{d}\rceil\geq\lceil\frac{pl+u_{b-i}-(n-1)}{d}\rceil$$
with equality holds if and only if
$$\frac{\tau(l)}{d}\geq\frac{pl+u_{b-i}}{d}-\lceil\frac{pl+u_{b-i}-(n-1)}{d}\rceil,$$
it suffices to show for each $1\leq i\leq b,$
$$
\sum_{l=0}^{n-1}(\lceil\frac{(p-1)l+u_{b-i}}{d}\rceil-\delta_{\in}^{(i)}(l))
=\sum_{l=0}^{n-1}\lceil\frac{pl+u_{b-i}-(n-1)}{d}\rceil.
$$
It is trivial for $n=1.$ For $n\geq2,$ it need to show
$$\sum_{l=1}^{n-1}(\delta_{<}^{(i)}(l)-\delta_{\in}^{(i)}(l))=\sharp\{0\leq l\leq
n-1|\{\frac{n-1}{d}\}^{'}<\{\frac{pl+u_{b-i}}{d}\}^{'}\}-1.$$
By Lemma \ref{mainrelation},
in case of $d|u_{b-i},$  the lemma follows from
$$\delta_{[1,n-1]}^{(i)}=1_{\frac{1}{d}\leq\{\frac{u_{b-i}}{d}\}\leq\frac{n-1}{d}}=0,$$
and
$$
\sharp\{0\leq l\leq
n-1|\{\frac{n-1}{d}\}^{'}<\{\frac{pl+u_{b-i}}{d}\}^{'}\}$$
$$=1+\sharp\{1\leq l\leq
n-1|\{\frac{n-1}{d}\}<\{\frac{pl+u_{b-i}}{d}\}\}.$$ In case of
$d\nmid u_{b-i},$ we have
$$
\sharp\{0\leq l\leq
n-1|\{\frac{n-1}{d}\}^{'}<\{\frac{pl+u_{b-i}}{d}\}^{'}\}$$
$$=\sharp\{1\leq l\leq
n-1|\{\frac{n-1}{d}\}<\{\frac{pl+u_{b-i}}{d}\}\}+\delta_{[1,n-1]}^{(i)}+1_{\frac{n-1}{d}<\{\frac{u_{b-i}}{d}\}}.
$$
Hence it follows from the equation
$1_{\frac{1}{d}\leq\{\frac{u_{b-i}}{d}\}\leq\frac{n-1}{d}}+1_{\frac{n-1}{d}<\{\frac{u_{b-i}}{d}\}}=1.$
\endproof

\begin{definition} Define $H_{n,u}(y)=\prod_{i=1}^bH_{n,u}^{(i)}(y),$ where
$$H_{n,u}^{(i)}(y)=\sum_{\tau\in S_{n,i}}{\rm
sgn}(\tau)\prod_{l\in A_n}
\sum\limits_{\stackrel{\sum\limits_{j=1}^djn_j=pl-\tau(l)+u_{b-i}}{\sum\limits_{j=1}^dn_j
=\lceil\frac{pl-\tau(l)+u_{b-i}}{d}\rceil}}
\prod_{j=1}^d\lambda_{n_j}y_j^{n_j}.$$
\end{definition}

\begin{theorem}\label{special-points}
Let $p>4d-\varepsilon(u),$ then we have
$$\prod_{i=1}^b\det((\gamma_{pl-j+u_{b-i}})_{l,j\in A_n})$$
$$=H_{n,u}((\hat{a}_i)_{i=0,1,\cdots,d})\pi^{bP_{[0,d],u}(n)}+O(\pi^{bP_{[0,d],u}(n)+\frac{1}{d(q-1)}}).$$
\end{theorem}
\proof
We have $$\det(\gamma_{pl-j+u_{b-i}})_{0\leq l,j\leq
n-1}=\sum_{\tau\in S_n}{\rm
sgn}(\tau)\prod_{l=0}^{n-1}\gamma_{pl-\tau(l)+u_{b-i}}.$$
The theorem now follows from Lemma \ref{rawesti} and \ref{estiofdet}.
\endproof

\begin{definition}
The reduction of $H_{n,u}$ modulo $p$ is denoted as
$\overline{H}_{n,u}$, and is called the $u$-twisted Hasse polynomial
of $[0,d]$ at $n$. The $u$-twisted Hasse polynomial $H_u$ of $[0,d]$
is defined by $H_u=\prod_{n=1}^{d-1}\overline{H}_{n,u}.$
\end{definition}
\begin{theorem}
The $u$-twisted Hasse polynomial $H_u$ of $[0,d]$ is nonzero.
\end{theorem}
\proof As a polynomial
of $y_d,$ the leading terms of $H_{n,u}^{(i)}(y)$ appear in
$$\sum_{\tau\in S_{n,i}}{\rm sgn}(\tau)\prod_{l=0}^{n-1}\lambda_{[\frac{pl+u_{b-i}-\tau(l)}{d}]}
y_d^{[\frac{pl+u_{b-i}-\tau(l)}{d}]}y_{d\{\frac{pl+u_{b-i}-\tau(l)}{d}\}}^{\delta_{\tau}(l)},$$
where $\delta_{\tau}(l)$ is 0 or 1 depending on whether
$d|(pl+u_{b-i}-\tau(l))$ or not.

We show among the leading
terms of $H_{n,u}^{(i)}(y),$ there is exactly one minimal monomial.

For $0\leq l\leq n-1,$ write $pl+u_{b-i}=q_ld+r_l,~~0\leq r_l\leq
d-1.$ Then $\tau\in S_{n,i}$ is equivalent to
$$\frac{\tau(l)}{d}\geq\frac{r_l}{d}-\lceil\frac{r_l-(n-1)}{d}\rceil.$$
Assume we have $r_{l_0}<r_{l_1}<\cdots<r_{l_m}\leq n-1<
r_{l_{m+1}}<\cdots<r_{l_{n-1}}.$ Hence $\tau\in S_{n,i}$ if and
only if $\tau(l_j)\geq r_{l_j}$ for $0\leq j\leq m.$

 Since
$$\sum_{l=0}^{n-1}[\frac{pl+u_{b-i}-\tau(l)}{d}]=\sum_{l=0}^{n-1}q_l+\sum_{j=0}^m[\frac{r_{l_j}-\tau(l_{j})}{d}]\leq\sum_{l=0}^{n-1}q_l,$$
with equality holding if and only if $r_{l_j}=\tau(l_{j})$ for all
$0\leq j\leq m.$ Therefore the leading terms appear for such $\tau$ that
$r_{l_j}=\tau(l_{j})$ for all $0\leq j\leq m.$ For $m+1\leq j\leq
n-1,$ we have $\delta_{\tau}(l_j)=1$ and
$\{\frac{pl_j+u_{b-i}-\tau(l_j)}{d}\}=\frac{r_{l_j}-\tau(l_j)}{d}.$
Hence among the leading terms, the minimal monomial appears exactly
when $\tau(l_{j})=r_{l_j}$ for all $0\leq j\leq m$ and

$$\tau(l_{m+1})=\max\{A_n-\{r_{l_0},\cdots,r_{l_m}\}\},$$
  $$\tau(l_{m+2})=\max\{A_n-\{r_{l_0},\cdots,r_{l_m},\tau(l_{m+1})\}\}, $$
  $$\cdots $$
  $$\tau(l_{n-1})=A_n-\{r_{l_0},\cdots,r_{l_m},\tau(l_{m+1}),\cdots,\tau(l_{n-2})\}.$$
Now the theorem follows from
$\lambda_j=\frac{1}{j!}\in\mathbb{Z}_p^{\times}$ for $0\leq
j\leq p-1.$
\endproof

\section{Proof of the main theorems}
In this section we prove the main theorems of this paper.

Firstly we prove Theorem \ref{main}, which says that $$P_{[0,d],u}\geq (p-1)H_{[0,d],u}^{\infty}$$
with equality holding at the point $d.$

{\it Proof of Theorem \ref{main}.} It need only to show this for $n\leq d.$
By the equation $$\sum\limits_{l=0}^{n-1}(\lceil\frac{(p-1)l+u_{i}}{d}\rceil-\delta_{\in}^{(i)}(l))=\sum\limits_{l=0}^{n-1}\lceil\frac{pl+u_{i}-(n-1)}{d}\rceil$$
that used in the proof of  Lemma \ref{estiofdet}, it suffices to show
$$\sum\limits_{l=0}^{n-1}\lceil\frac{pl+u_{i}-(n-1)}{d}\rceil\geq\sum_{l=0}^{n-1}\frac{(p-1)l+u_{i}}{d}.$$ The case $n=1$  is trivial. For $n>1,$
since $$
\sum\limits_{l=0}^{n-1}\lceil\frac{pl+u_{i}-(n-1)}{d}\rceil
=\sum_{l=0}^{n-1}(\frac{pl+u_i}{d}-\{\frac{pl+u_i}{d}\}^{'}+1_{\{\frac{n-1}{d}\}<\{\frac{pl+u_i}{d}\}^{'}}),
$$
it suffices to show
$$\sum_{l=0}^{n-1}(\{\frac{pl+u_i}{d}\}^{'}-1_{\{\frac{n-1}{d}\}<\{\frac{pl+u_i}{d}\}^{'}})\leq\sum_{l=0}^{n-1}\frac{l}{d}.$$
Now the inequality follows from
$$\sum_{l=0}^{n-1}(\{\frac{pl+u_i}{d}\}^{'}-1_{\{\frac{n-1}{d}\}<\{\frac{pl+u_i}{d}\}^{'}})
\leq\sum_{l=0,\{\frac{pl+u_i}{d}\}^{'}\leq\{\frac{n-1}{d}\}}^{n-1}\{\frac{pl+u_i}{d}\}^{'}.$$
Moreover the equalities above hold when $n=d,$
the theorem is proved.\qed

\begin{lemma}\label{zqtozp}
The $T$-adic Newton polygon of
$\det(1-\Psi^as^a|B/\mathbb{Z}_q[[\pi^{\frac{1}{d(q-1)}}]])$
coincides with that of $\det(1-\Psi
s|B/\mathbb{Z}_p[[\pi^{\frac{1}{d(q-1)}}]]).$
\end{lemma}

\proof The lemma follows from the following:
\begin{align*}
\prod_{\zeta^a=1}\det(1-\Psi\zeta
s|B/\mathbb{Z}_p[[\pi^{\frac{1}{d(q-1)}}]])&=\det(1-\Psi^as^a|B/\mathbb{Z}_p[[\pi^{\frac{1}{d(q-1)}}]])\\
&={\rm
Norm}(\det(1-\Psi^as^a|B/\mathbb{Z}_q[[\pi^{\frac{1}{d(q-1)}}]])),
\end{align*}
where Norm is the norm map from
$\mathbb{Q}_q[[\pi^{\frac{1}{d(q-1)}}]]$ to
$\mathbb{Q}_p[[\pi^{\frac{1}{d(q-1)}}]].$
\endproof

\begin{lemma}\label{BtoB_i}
The T-adic Newton polygon of $C_{f,u}(s,T)^b$ coincides with that of
$\det(1-\Psi^as|B/\mathbb{Z}_q[[\pi^{\frac{1}{d(q-1)}}]]).$
\end{lemma}

\proof Let $\text{Gal}(\mathbb{Q}_q/\mathbb{Q}_p)$ act on $\mathbb{Z}_q[[T]][[s]]$ by fixing $s$ and $T,$
then we have $C_{pu,f}(s,T)=C_{f,u}(s,T)^{\sigma}.$
Therefore the lemma follows from the following
$$\prod_{i=0}^{b-1}C_{f,u}(s,T)^{\sigma^i}=\prod_{i=0}^{b-1}C_{p^iu,f}(s,T)$$
$$=\prod_{i=0}^{b-1}\det(1-\Psi^as|B_{p^iu}/\mathbb{Z}_q[[\pi^{\frac{1}{d(q-1)}}]])
=\det(1-\Psi^as|B/\mathbb{Z}_q[[\pi^{\frac{1}{d(q-1)}}]]).$$
\endproof

\begin{theorem}\label{description of NP}
The $T$-adic Newton polygon of $C_{f,u}(s,T)$ is the lower convex closure of the points
$$(n,\frac{1}{b}\text{ord}_{T}(C_{abn})),~n=0,1,\cdots.$$
\end{theorem}
\proof By Lemma \ref{zqtozp}, we see the T-adic Newton polygon of
the power series
$\det(1-\Psi^as^a|B/\mathbb{Z}_q[[\pi^{\frac{1}{d(q-1)}}]])$ is the
lower convex closure of the points $$(n,{\rm
ord}_{T}(C_{n})),~~~n=0,1,\cdots.$$ It is clear that
$(n,\text{ord}_{T}(C_n))$ is not a vertex of that polygon if $a\nmid
n$. So that Newton polygon is the lower convex closure of the points
$$(an,\text{ord}_{T}(C_{an})),\ n=0,1,\cdots.$$
Hence the T-adic Newton polygon of
$\det(1-\Psi^as|B/\mathbb{Z}_q[[\pi^{\frac{1}{d(q-1)}}]])$ is the
convex closure of the points $$(n,{\rm
ord}_{T}(C_{an})),~~~n=0,1,\cdots.$$ By Lemma \ref{BtoB_i}, the T-adic
Newton polygon of $C_{f,u}(s,T)^b$ is the lower convex closure of
the points $$(n,{\rm ord}_{T}(C_{an})),~~~n=0,1,\cdots,$$ hence the
closure of the points $$(bn,{\rm ord}_{T}(C_{abn})),~~~n=0,1,\cdots.$$
It follows that the $T$-adic Newton polygon of $C_{f,u}(s,T)$ is the
convex closure of the points
$$(n,\frac{1}{b}\text{ord}_{T}(C_{abn})),~n=0,1,\cdots.$$ The theorem is proved.
\endproof
We now prove Theorem \ref{main1}, which says that if $p>4d-\varepsilon(u),$
$$T-\text{adic NP of }C_{f,u}(s,T)\geq\text{ord}_p(q)P_{[0,d],u}. $$
{\it Proof of Theorem \ref{main1}.} The theorem follows from
Theorem \ref{order of C_bn} and the last theorem.\qed

\begin{theorem}\label{spe} Let $A(s,T)$ be a
$T$-adic entrie series in $s$ with unitary constant term. If
$0\neq|t|_p<1$, then
$$t-adic\text{ NP of
}A(s,t)\geq T-adic\text{ NP of }A(s,T),$$ where NP is the short for
Newton polygon. Moreover, the equality holds for one $t$ if and only
if it holds for all $t$.\end{theorem}
\proof
The reader may refer \cite{LLN} and we omit the proof here.
\endproof

\begin{theorem}Let $f(x)=\sum\limits_{i=0}^d(a_ix^i,0,0,\cdots)$, and $p>4d-\varepsilon(u).$
If the equality $$\pi_m-\text{adic NP of }C_{f,u}(s,\pi_m)={\rm
ord}_p(q)P_{[0,d],u}$$holds for one $m\geq1$, then it holds for all
$m\geq1$, and we have $$T-\text{adic NP of }C_{f,u}(s,T)={\rm
ord}_p(q)P_{[0,d],u}.$$
\end{theorem}
\proof
This follows from Theorem \ref{main1} and the last theorem.
\endproof

\begin{theorem}Let $f(x)=\sum\limits_{i=0}^d(a_ix^i,0,0,\cdots)$.
Then  $$\pi_m-\text{adic NP of }C_{f,u}(s,\pi_m)={\rm
ord}_p(q)P_{[0,d],u}$$ if and only if
$$\pi_m-\text{adic NP of }L_{f,u}(s,\pi_m)={\rm
ord}_p(q)P_{[0,d],u}\text{ on }[0,p^{m-1}d].$$
\end{theorem}
\proof Assume that $ L_{f,u}(s,
\pi_m)=\prod\limits_{i=1}^{p^{m-1}d}(1-\beta_is) $. Then
$$C_{f,u}(s,\pi_m )= \prod\limits_{j=0}^{\infty} L_{f,u}(q^js, \pi_m)=
\prod\limits_{j=0}^{\infty}\prod\limits_{i=1}^{p^{m-1}d}(1-\beta_iq^js).$$
Therefore the slopes of the  $q$-adic Newton polygon of
$C_{f,u}(s,\pi_m)$ are the numbers
$$j+{\rm ord}_q(\beta_i),\ 1\leq
i\leq p^{m-1}{\rm Vol}(\triangle), j=0,1,\cdots.$$
Since $$w(i+m^{m-1}d)=p^m-p^{m-1}+w(i),$$
then the slopes of $P_{[0,d],u}$ are the numbers
 $$j(p^m-p^{m-1})+w(i),\ 1\leq i\leq
p^{m-1}d, j=0,1,\cdots.$$ It follows
that
$$\pi_m-\text{adic NP of }C_{f,u}(s,\pi_m)={\rm
ord}_p(q)P_{[0,d],u}$$ if and only if
$$\pi_m-\text{adic NP of }L_{f,u}(s,\pi_m)={\rm
ord}_p(q)P_{[0,d],u}\text{ on }[0,p^{m-1}d].$$
\endproof

We now prove Theorems \ref{main2}, \ref{main3} and \ref{main4}. By
the above theorems, it suffices to prove the following.
\begin{theorem}\label{final th}
Let $f(x)=\sum\limits_{i=0}^d(a_ix^i,0,0,\cdots)$, and $p>4d-\varepsilon(u)$.
Then  $$\pi_1-\text{adic NP of }L_{f,u}(s,\pi_1)={\rm
ord}_p(q)P_{[0,d],u}\text{ on }[0,d]$$ if and
only if $H((a_i)_{0\leq i\leq d})\neq0$.
\end{theorem}\proof By a result of Liu\cite{Liu2}, the
$q$-adic Newton polygon  of $L_{f,u}(s,\pi_1)$ coincides with
$H_{[0,d],u}^{\infty}$ at the point $d$. By
Theorem \ref{main}, $P_{[0,d],u}$ coincides with
$(p-1)H_{[0,d],u}^{\infty}$ at the point $d$,
It follows that the $\pi_1$-adic Newton polygon of $L_{f,u}(s,\pi_1)$
coincides with ${\rm ord}_p(q)P_{[0,d],u}$ at the point $d$. Therefore it suffices to show that
$$\pi_1-\text{adic NP of }L_{f,u}(s,\pi_1)={\rm
ord}_p(q)P_{[0,d],u}\text{ on }[0,d-1]$$ if and
only if $H((a_i)_{0\leq i\leq d})\neq0$.

From the identity $$C_{f,u}(s,\pi_1 )= \prod\limits_{j=0}^{\infty}
L_{f,u}(q^js, \pi_1),$$ and the fact the $q$-adic orders of the
reciprocal roots of $L_{f,u}(s,\pi_1)$ are no greater than $1$, we
infer that
$$\pi_1-\text{adic NP of }L_{f,u}(s,\pi_1) =\pi_1-\text{adic NP of
}C_{f,u}(s,\pi_1)\text{ on }[0,d-1].$$Therefore it
suffices to show that $$\pi_1-\text{adic NP of }C_{f,u}(s,\pi_1)={\rm
ord}_p(q)P_{[0,d],u}\text{ on }[0,d-1]$$ if and
only if $H((a_i)_{0\leq i\leq d})\neq0$. The theorem now follows
from the $T$-adic Dwork trace formula and Theorems
\ref{order of C_bn} and \ref{special-points}.\endproof



\begin{thebibliography}{99}
\bibitem{AS1}A. Adolphson and S. Sperber,
Exponential sums and Newton polyhedra: cohomology and estimates, Ann. Math., 130 (1989), 367-406.

\bibitem{AS2}A. Adolphson and S. Sperber,
Newton polyhedra and the degree of the L-function associated to an exponential sum, Invent. Math. 88(1987), 555-569.

\bibitem{AS3}A. Adolphson and S. Sperber,
On twisted exponential sums, Math. Ann., 290 (1991), 713-726.

\bibitem{AS4}A. Adolphson and S. Sperber, Twisted exponential sums and
Newton polyhedra, J. reine angew. Math., 443 (1993), 151-177.

\bibitem{BF}R. Blache and E. F\'erard, Newton straitification for
polynomials: the open stratum, J. Number Theory, 123(2007), 456-472.

\bibitem{BFZ}R. Blache, E. F\'erard and J.H. Zhu,
Hodge-Stickelberger polygons for L-functions of exponential sums of
$P(x^s)$, Math. Res. Letter, 15 (2008), 1053-1071.

\bibitem{Lw}W. -C. W. Li, Character sums over finite fields, J. Number Theory 74(1999), 181-229.

\bibitem{LZ}Hanfeng Li and Hui June Zhu, Zeta functions of totally
ramified p-covers of the projective line, Rend. Sem. Mat. Univ.
Padova, Vol. 113(2005), 203-225.

\bibitem{Liu1}C.Liu, The L-functions of twisted Witt coverings, J. Number Theory, 125(2007), 267-284.

\bibitem{Liu2}C.Liu, T-adic exponential sums under base change, arxiv:0909.1111.

\bibitem{LWn}C. Liu and D. Wan, $T$-adic exponential sums, Algebra \& Number Theory, Vol. 3, No. 5 (2009), 489-509.

\bibitem{LW}C. Liu and D. Wei, The $L$-functions of Witt
coverings, Math. Z., 255 (2007), 95-115.

\bibitem{LLN}C. Liu, Wenxin Liu and Chuanze Niu, T-adic exponential
sums in one variable, arxiv:0901.0354v5.

\bibitem{Se}J-P. Serre, Endomorphismes compl\'etement continus des
espaces de Banach $p$-adiques, Publ. Math., IHES., 12(1962), 69-85.

\bibitem{DWn} D. Wan, Variation of $p$-adic Newton polygons for L-functions of exponential
sums, Asian J. Math., Vol 8, 3(2004), 427-474.
\bibitem{Zh1}J. H. Zhu, p-adic variation of L functions of one variable exponential sums,
I. Amer. J.  Math., 125 (2003), 669-690.

\bibitem{Zh2}J. H. Zhu, Asymptotic variation of L functions of one-variable exponential
sums, J. Reine Angew. Math., 572 (2004), 219--233.

\bibitem{Zhu3}J. H. Zhu, L-functions of exponential sums over one-dimensional affinoids :
Newton over Hodge, Inter. Math. Research Notices, no 30 (2004),
1529--1550.

\end{thebibliography}
\end{document}